\documentclass[11pt]{amsart}
\usepackage{pstricks,amssymb}

\def \blanc{ \hskip1mm}

\newcommand{\sa}{\begin{pspicture}[.2](0,0)(.15,.3)
\psline[linewidth=.004]{-}(0,.3)(.3,.3)
\psline[linewidth=.004]{-}(0,.2)(.3,.2)
\psline[linewidth=.004]{-}(0,.1)(.3,.0)
\psline[linewidth=.004]{-}(0,.0)(.1,.015)
\psline[linewidth=.004]{-}(.2,.085)(.3,.1)
\end{pspicture}
}

\newcommand{\saa}{\begin{pspicture}[.2](0,0)(.15,.3)
\psline[linewidth=.004]{-}(0,.3)(.3,.3)
\psline[linewidth=.004]{-}(0,.2)(.3,.2)
\psline[linewidth=.004]{-}(0,.1)(.1,.085)
\psline[linewidth=.004]{-}(0,.0)(.3,.1)
\psline[linewidth=.004]{-}(.2,.015)(.3,.0)
\end{pspicture}
}

\newcommand{\sbbb}{\begin{pspicture}[.2](0,0)(.15,.3)
\psline[linewidth=.004]{-}(0,.3)(.3,.3)
\psline[linewidth=.004]{-}(0,0)(.3,.0)
\psline[linewidth=.004]{-}(0,.2)(.3,.1)
\psline[linewidth=.004]{-}(0,.1)(.1,.115)
\psline[linewidth=.004]{-}(.2,.185)(.3,.2)
\end{pspicture}
}

\newcommand{\sbb}{\begin{pspicture}[.2](0,0)(.15,.3)
\psline[linewidth=.004]{-}(0,.3)(.3,.3)
\psline[linewidth=.004]{-}(0,0)(.3,.0)
\psline[linewidth=.004]{-}(0,.2)(.1,.185)
\psline[linewidth=.004]{-}(0,.1)(.3,.2)
\psline[linewidth=.004]{-}(.2,.115)(.3,.1)
\end{pspicture}
}

\newcommand{\sccc}{\begin{pspicture}[.2](0,0)(.15,.3)
\psline[linewidth=.004]{-}(0,.1)(.3,.1)
\psline[linewidth=.004]{-}(0,0)(.3,.0)
\psline[linewidth=.004]{-}(0,.2)(.1,.215)
\psline[linewidth=.004]{-}(0,.3)(.3,.2)
\psline[linewidth=.004]{-}(.2,.285)(.3,.3)
\end{pspicture}
}

\newcommand{\scc}{\begin{pspicture}[.2](0,0)(.15,.3)
\psline[linewidth=.004]{-}(0,.1)(.3,.1)
\psline[linewidth=.004]{-}(0,0)(.3,.0)
\psline[linewidth=.004]{-}(0,.2)(.3,.3)
\psline[linewidth=.004]{-}(0,.3)(.1,.285)
\psline[linewidth=.004]{-}(.2,.215)(.3,.2)
\end{pspicture}
}

\newcommand{\Saaa}{\begin{pspicture}[.2](0,0)(.15,.5)
\psline[linewidth=.004]{-}(0,.4)(.3,.4)
\psline[linewidth=.004]{-}(0,.3)(.3,.3)
\psline[linewidth=.004]{-}(0,.2)(.3,.2)
\psline[linewidth=.004]{-}(0,.1)(.3,.0)
\psline[linewidth=.004]{-}(0,.0)(.1,.015)
\psline[linewidth=.004]{-}(.2,.085)(.3,.1)
\end{pspicture}
}

\newcommand{\Saa}{\begin{pspicture}[.2](0,0)(.15,.5)
\psline[linewidth=.004]{-}(0,.4)(.3,.4)
\psline[linewidth=.004]{-}(0,.3)(.3,.3)
\psline[linewidth=.004]{-}(0,.2)(.3,.2)
\psline[linewidth=.004]{-}(0,.1)(.1,.085)
\psline[linewidth=.004]{-}(0,.0)(.3,.1)
\psline[linewidth=.004]{-}(.2,.015)(.3,.0)
\end{pspicture}
}

\newcommand{\Sbbb}{\begin{pspicture}[.2](0,0)(.15,.5)
\psline[linewidth=.004]{-}(0,.4)(.3,.4)
\psline[linewidth=.004]{-}(0,.3)(.3,.3)
\psline[linewidth=.004]{-}(0,0)(.3,.0)
\psline[linewidth=.004]{-}(0,.2)(.3,.1)
\psline[linewidth=.004]{-}(0,.1)(.1,.115)
\psline[linewidth=.004]{-}(.2,.185)(.3,.2)
\end{pspicture}
}

\newcommand{\Sbb}{\begin{pspicture}[.2](0,0)(.15,.5)
\psline[linewidth=.004]{-}(0,.4)(.3,.4)
\psline[linewidth=.004]{-}(0,.3)(.3,.3)
\psline[linewidth=.004]{-}(0,0)(.3,.0)
\psline[linewidth=.004]{-}(0,.2)(.1,.185)
\psline[linewidth=.004]{-}(0,.1)(.3,.2)
\psline[linewidth=.004]{-}(.2,.115)(.3,.1)
\end{pspicture}
}

\newcommand{\Sccc}{\begin{pspicture}[.2](0,0)(.15,.5)
\psline[linewidth=.004]{-}(0,.4)(.3,.4)
\psline[linewidth=.004]{-}(0,.1)(.3,.1)
\psline[linewidth=.004]{-}(0,0)(.3,.0)
\psline[linewidth=.004]{-}(0,.2)(.1,.215)
\psline[linewidth=.004]{-}(0,.3)(.3,.2)
\psline[linewidth=.004]{-}(.2,.285)(.3,.3)
\end{pspicture}
}

\newcommand{\Scc}{\begin{pspicture}[.2](0,0)(.15,.5)
\psline[linewidth=.004]{-}(0,.4)(.3,.4)
\psline[linewidth=.004]{-}(0,.1)(.3,.1)
\psline[linewidth=.004]{-}(0,0)(.3,.0)
\psline[linewidth=.004]{-}(0,.2)(.3,.3)
\psline[linewidth=.004]{-}(0,.3)(.1,.285)
\psline[linewidth=.004]{-}(.2,.215)(.3,.2)
\end{pspicture}
}

\newcommand{\Sddd}{\begin{pspicture}[.2](0,0)(.15,.5)
\psline[linewidth=.004]{-}(0,0)(.3,0)
\psline[linewidth=.004]{-}(0,.2)(.3,.2)
\psline[linewidth=.004]{-}(0,.1)(.3,.1)
\psline[linewidth=.004]{-}(0,.3)(.1,.315)
\psline[linewidth=.004]{-}(0,.4)(.3,.3)
\psline[linewidth=.004]{-}(.2,.385)(.3,.4)
\end{pspicture}
}

\newcommand{\Sdd}{\begin{pspicture}[.2](0,0)(.15,.5)
\psline[linewidth=.004]{-}(0,0)(.3,0)
\psline[linewidth=.004]{-}(0,.2)(.3,.2)
\psline[linewidth=.004]{-}(0,.1)(.3,.1)
\psline[linewidth=.004]{-}(0,.3)(.3,.4)
\psline[linewidth=.004]{-}(0,.4)(.1,.385)
\psline[linewidth=.004]{-}(.2,.315)(.3,.3)
\end{pspicture}
}

\newcommand{\fcp}{\begin{pspicture}[.2](0,0)(.5,.4)
\psline{-}(0.05,0)(.45,.4)
\psline{-}(.45,0)(.3,.15)
\psline{-}(.2,.25)(0.05,.4)
\end{pspicture}
}

\newcommand{\fcm}{\begin{pspicture}[.2](0,0)(.5,.4)
\psline{-}(0.45,0)(.05,.4)
\psline{-}(.05,0)(.2,.15)
\psline{-}(.3,.25)(0.45,.4)
\end{pspicture}
}

\newcommand{\fcpo}{\begin{pspicture}[.2](0,0)(.5,.4)
\psline{-}(0.05,0)(.45,.4)
\psline{-}(.45,0)(.3,.15)
\psline{-}(.2,.25)(0.05,.4)
\psline[linewidth=.004]{-}(0.05,.25)(0.05,.4)
\psline[linewidth=.004]{-}(0.20,.4)(0.05,.4)
\psline[linewidth=.004]{-}(0.30,0.4)(0.45,0.4)
\psline[linewidth=.004]{-}(0.45,0.25)(0.45,0.4) 
\end{pspicture}
}

\newcommand{\fcmo}{\begin{pspicture}[.2](0,0)(.5,.4)
\psline{-}(0.45,0)(.05,.4)
\psline{-}(.05,0)(.2,.15)
\psline{-}(.3,.25)(0.45,.4)
\psline[linewidth=.004]{-}(0.05,.25)(0.05,.4)
\psline[linewidth=.004]{-}(0.20,.4)(0.05,.4)
\psline[linewidth=.004]{-}(0.30,0.4)(0.45,0.4)
\psline[linewidth=.004]{-}(0.45,0.25)(0.45,0.4) 
\end{pspicture}
}

\newcommand{\postref}{\begin{pspicture}[.2](-.1,0)(1.3,.6)
\pscurve{-}(1.05,.35)(1.2,.4)(1,.5)(.6,.6)(.2,.5)(0,.4)(.2,.3)(.4,.15)(.55,.25)
\pscurve{-}(.25,.35)(.4,.45)(.6,.3)(.8,.15)(.95,.25)
\pscurve{-}(.15,.25)(.0,.2)(.2,.1)(.6,0)(1,.1)(1.2,.2)(1,.3)(.8,.45)(.65,.35)
\end{pspicture}}

\newcommand{\cuspgp}{\begin{pspicture}[.2](-0.1,0)(.6,.4)
\pscurve{-}(.0,.2)(.1,.19)(.2,.17)(.3,.14)(.4,.10)(.5,.0)
\pscurve{-}(.0,.2)(.1,.21)(.2,.23)(.3,.26)(.4,.30)(.5,.4)
\psline[linewidth=.004]{-}(.5,.25)(.5,.4)
\psline[linewidth=.004]{-}(.35,.4)(.5,.4)
\end{pspicture}}

\newcommand{\cuspg}{\begin{pspicture}[.2](-0.1,0)(.6,.4)
\pscurve{-}(.0,.2)(.1,.19)(.2,.17)(.3,.14)(.4,.10)(.5,.0)
\pscurve{-}(.0,.2)(.1,.21)(.2,.23)(.3,.26)(.4,.30)(.5,.4)
\end{pspicture}}

\newcommand{\cuspd}{\begin{pspicture}[.2](-0.1,0)(.7,.4)
\pscurve{-}(.5,.2)(.4,.19)(.3,.17)(.2,.14)(.1,.10)(.0,.0)
\pscurve{-}(.5,.2)(.4,.21)(.3,.23)(.2,.26)(.1,.30)(.0,.4)

\end{pspicture}}

\newcommand{\cuspdn}{\begin{pspicture}[.2](-0.1,0)(.7,.4)
\pscurve{-}(.5,.2)(.4,.19)(.3,.17)(.2,.14)(.1,.10)(.0,.0)
\pscurve{-}(.5,.2)(.4,.21)(.3,.23)(.2,.26)(.1,.30)(.0,.4)
\psline[linewidth=.004]{-}(.0,.15)(.0,0)
\psline[linewidth=.004]{-}(.15,0)(.0,.0)
\end{pspicture}}

\newcommand{\boucled}{\begin{pspicture}[.2](-.1,0)(.9,.4)
\psline{-}(0,0)(.3,.3)
\pscurve{-}(.3,.3)(.5,.4)(.8,.2)(.5,0)(.3,.1)
\psline{-}(.3,.1)(.25,.15)
\psline{-}(.15,.25)(0,.4)
\end{pspicture}}

\newcommand{\droit}{\begin{pspicture}[.2](-0.1,0)(.5,.4)
\pscurve{-}(.0,.4)(.1,.38)(.2,.35)(.3,.29)(.33,.2)(.3,.11)(.2,.05)(.1,.02)(.0,0)
\end{pspicture}}

\newcommand{\gauche}{\begin{pspicture}[.2](-0.1,0)(.5,.4)
\pscurve{-}(.4,.4)(.3,.38)(.2,.35)(.1,.29)(.07,.2)(.1,.11)(.2,.05)(.3,.02)(.4,0)
\end{pspicture}}

\newcommand{\rond}{\begin{pspicture}[.2](-0.1,0)(.5,.4)
\pscircle(.2,.2){.2}
\end{pspicture}}

\newcommand{\ouvertV}{\begin{pspicture}[.2](-0.1,0)(.5,.4)
\pscurve{-}(0,0)(.1,.1)(.15,.2)(.1,.3)(0,.4)
\pscurve{-}(.4,.4)(.3,.3)(.25,.2)(.3,.1)(.4,0)
\end{pspicture}}

\newcommand{\ouvertVff}{\begin{pspicture}[.2](-0.1,0)(.5,.4)
\pscurve{-}(0,0)(.1,.1)(.15,.2)(.1,.3)(0,.4)
\pscurve{-}(.4,.4)(.3,.3)(.25,.2)(.3,.1)(.4,0)
\psline[linewidth=.004]{-}(0,.25)(0,.4)
\psline[linewidth=.004]{-}(0.15,.4)(0,.4)
\psline[linewidth=.004]{-}(0.25,0.4)(0.4,0.4)
\psline[linewidth=.004]{-}(0.4,0.25)(0.4,0.4) 
\end{pspicture}}

\newcommand{\ouvertVn}{\begin{pspicture}[.2](-0.1,0)(.5,.4)
\pscurve{-}(0,0)(.1,.1)(.15,.2)(.1,.3)(0,.4)
\pscurve{-}(.4,.4)(.3,.3)(.25,.2)(.3,.1)(.4,0)
\psline[linewidth=.004]{-}(0,.15)(0,0)
\psline[linewidth=.004]{-}(0.15,0)(0,0)
\psline[linewidth=.004]{-}(0.25,0.4)(0.4,0.4)
\psline[linewidth=.004]{-}(0.4,0.25)(0.4,0.4) 
\end{pspicture}}

\newcommand{\ouvertH}{\begin{pspicture}[.2](-0.1,0)(.5,.4)
\pscurve{-}(0,0)(.1,.1)(.2,.15)(.3,.1)(.4,0)
\pscurve{-}(.4,.4)(.3,.3)(.2,.25)(.1,.3)(0,.4)
\end{pspicture}}

\newcommand{\ouvertHp}{\begin{pspicture}[.2](-0.1,0)(.5,.4)
\pscurve{-}(0,0)(.1,.1)(.2,.15)(.3,.1)(.4,0)
\pscurve{-}(.4,.4)(.3,.3)(.2,.25)(.1,.3)(0,.4)
\psline[linewidth=.004]{-}(0,.15)(0,0)
\psline[linewidth=.004]{-}(.15,0)(0,0)
\psline[linewidth=.004]{-}(0.25,0.4)(0.4,0.4)
\psline[linewidth=.004]{-}(0.4,0.25)(0.4,0.4) 
\end{pspicture}}

\newcommand{\cross}{\begin{pspicture}[.2](-0.1,0)(.5,.5)
\psline{-}(0,0)(.4,.4)
\psline{-}(0,.4)(.4,0)
\end{pspicture}}

\newcommand{\infi}{\begin{pspicture}[.2](-0.1,0)(.9,.4)
\pscurve{-}(.45,.15)(.6,0)(.8,.2)(.6,.4)(.4,.2)(.2,0)(0,.2)(.2,.4)(.35,.25)
\end{pspicture}}

\newcommand{\trifi}{\begin{pspicture}[.2](-0.1,0)(1.3,.4)
\pscurve{-}(.35,.25)(.2,.4)(0,.2)(.2,0)(.4,.2)(.6,.4)(.75,.25)
\pscurve{-}(.45,.15)(.6,.0)(.8,.2)(1,.4)(1.2,.2)(1,0)(.85,.15)
\end{pspicture}}

\newcommand{\flying}{\begin{pspicture}[.2](-0.1,0)(.7,.4)
\pscurve{-}(0,.2)(.2,.25)(.3,.35)(.4,.25)(.6,.2)
\pscurve{-}(0,.2)(.2,.15)(.3,.05)(.4,.15)(.6,.2)
\end{pspicture}}

\newcommand{\bifly}{\begin{pspicture}[.2](-0.1,0)(1.1,.4)
\pscurve{-}(0,.2)(.15,.25)(.3,.35)(.5,.2)(.7,.05)(.85,.15)(1,.2)
\pscurve{-}(0,.2)(.15,.15)(.3,.05)(.5,.2)(.7,.35)(.85,.25)(1,.2)
\end{pspicture}}

\newcommand{\statun}{\begin{pspicture}[.2](-0.1,0)(1.1,.4)
\pscurve{-}(0,.2)(.15,.25)(.3,.4)(.5,.2)(.7,.0)(.85,.15)(1,.2)
\pscurve{->}(0,.2)(.15,.15)(.3,.0)(.4,.09)(.45,.15)
\pscurve{-}(.45,.15)(.5,.2)(.7,.4)(.85,.25)(1,.2)
\end{pspicture}}

\newcommand{\statdeux}{\begin{pspicture}[.2](-0.1,0)(1.1,.4)
\pscurve{->}(0,.2)(.15,.25)(.3,.4)(.4,.31)(.45,.25)
\pscurve{-}(.45,.25)(.5,.2)(.7,.0)(.85,.15)(1,.2)
\pscurve{-}(0,.2)(.15,.15)(.3,.0)(.5,.2)(.7,.4)(.85,.25)(1,.2)
\end{pspicture}}

\newcommand{\stattrois}{\begin{pspicture}[.2](-0.1,0)(1,.4)
\pscurve{->}(0,.2)(.2,.25)(.3,.32)(.5,.35)
\pscurve{-}(1,.2)(.8,.25)(.7,.32)(.5,.35)
\pscurve{-}(0,.2)(.2,.15)(.5,.05)(.8,.15)(1,.2)
\end{pspicture}}

\font\tmsb=msbm10 at12pt
\font\smsb=msbm7
\font\ssmsb=msbm5
\newfam\msbfam
\textfont\msbfam=\tmsb
\scriptfont\msbfam=\smsb
\scriptscriptfont\msbfam=\ssmsb
\def \mth{\fam\msbfam}
\def \Mth#1{{\mth #1}}
\def \RM{\Mth{R}}

\title{On Legendrian knots and polynomial invariants}

\author{Emmanuel Ferrand}
\address{Institut Fourier, BP 74, 38402 St Martin d'H\`eres Cedex, France.}

\email{emmanuel.ferrand@ujf-grenoble.fr}

\date{First version: February 1st, 2000; this 3rd version: July 10, 2000.}

\begin{document}

\begin{abstract}
It is proved in this note that the analogues of the Bennequin
inequality which provide an upper bound for
the Bennequin invariant of a Legendrian knot in the standard contact
three dimensional space in terms of the least degree in the framing
variable of the HOMFLY and the Kauffman polynomials are not sharp.
Furthermore, the relationships between these restrictions on the range
of the Bennequin invariant are investigated, which leads to a  
new simple proof of the inequality involving
the Kauffman polynomial.
\end{abstract}

\maketitle


\bigskip

\section{Introduction}
The {\em standard contact three dimensional space} is $\RM^3$ with
coordinates $u,p,q$ endowed with
the plane field induced by the the contact form $\alpha=du-pdq$ and
the orientation induced by the volume form $du \wedge dp \wedge dq$.
A smooth knot embedded in  $\RM^3$ is called {\em Legendrian} 
if it is everywhere tangent to this
plane field. A Legendrian knot is completely determined by its
projection to the plane $(p,q)$, the other coordinate $u$ being the
integral
of the form $pdq$ along the knot projection. Any smooth plane curve
unambiguously defines a Legendrian immersion in $\RM^3$ provided that the integral of
$pdq$ vanishes along this curve (so that the Legendrian ``lift'' is
closed).
Any Legendrian knot is {\em horizontal} (i.e. the tangent vector is never vertical) 
with respect to this projection.
Hence a {\em contact isotopy} (a one parameter family of Legendrian knots)
is a particular case of what is classically called a {\em regular
  isotopy} (the first Reidemeister move is forbidden on knot
projections).  
Two horizontal oriented knots are regular isotopic if an only if they are
isotopic and have the same {\em writhe} $w$ (self-linking number with
respect to the vertical framing) and the same {\em Whitney index} $r$
(degree of the Gauss map of the knot projection).
Let $l$ be some Legendrian knot.
By definition, its Bennequin invariant is 
$tb(l)=-w(l)$ and its Maslov invariant is $\mu (l)=r(l)$.

\bigskip

The following restrictions are known (see below for a list of
the authors of these results):

\bigskip
{\bf Theorem.}
\begin{itemize}
\item[a).] $tb(l)+|\mu(l)| \leq 2.g_4(l)-1$
\item[b).] $tb(l)+|\mu(l)| \leq e_P(l)$
\item[c).] $tb(l) \leq e_Y(l)$
\end{itemize}

\bigskip

Here $g_4(l)$ denotes the slice genus of $l$, $e_P(l)$ (resp. $e_Y(l)$)
the least degree of the framing variable $a$ in the HOMFLY
(resp. Kauffman) polynomial of $l$ (see below for the precise
definition and normalization of these polynomials).

\bigskip

In this paper, the relationships between these three
inequalities are investigated. A new proof of $c)$ is given. 
It is shown that in spite of the misleading evidence provided by the
knot tables, there is no inequality like $e_P \leq 2g_4-1$,
from which it would follow that $b)$ implies $a)$.
We provide examples showing that inequality $e_Y \leq e_P$ is false,
again in spite of what could be expected from the tables.
As a consequence, none of the three inequalities above is sharp.

\section{Known results.}
\label{known}

\subsection{Inequalities}

\begin{itemize}

\item Consider a braid $\sigma$ with $n(\sigma)$ strands and whose exponent sum is $c(\sigma)$.
Denote by $\hat{\sigma}$ the closure of this braid and let $I$ be some knot invariant 
that does not detect the orientation of knots.
It follows from \cite{Be} (theorem 8, proposition 6, and paragraph 24) that if the inequality\footnote{The sign difference with \cite{Be} is due to the fact that we use a different contact form and a different orientation.}
$-c(\sigma)-n(\sigma) \leq I(\hat{\sigma})$ holds for any braid $\sigma$, 
then, for any Legendrian knot $l$ having 
topological knot type $K$, $tb(l)+|\mu(l)| \leq I(K)$.

\bigskip

\item Bennequin \cite{Be} (theorem 3) proved that $|c(\sigma)|-n(\sigma) \leq 2g_3(\hat{\sigma})-1$, where $g_3$ 
denotes the genus. Hence, by the previous discussion, $tb(l)+|\mu(l)| \leq 2g_3(l)-1$ 
(\cite{Be}, theorem 11).

\bigskip

\item In \cite{Ru5}, Lee Rudolph proved that, for any braid $\sigma$,  $|c(\sigma)|-n(\sigma) \leq 2g_4(\hat{\sigma})-1$, 
and hence inequality a) follows.

\bigskip

\item In \cite{Mo} an \cite{FW}, Morton and Franks-Williams proved, using ``elementary'' combinatorial means,
that\footnote{We follow \cite{Kau} for the normalization of the Homfly
  polynomial. This explains the difference with the original 
inequality of \cite{Mo},  where different variables and a different normalization are assumed.}  
$-c(\sigma)-n(\sigma) \leq e_P(\hat{\sigma})$. 
As observed in \cite{FT}, inequality b) follows from this and the preceding discussion.

\bigskip

\item An inequality similar to inequality c) with $e_Y$
replaced by the least degree of the framing variable in the Kauffman
polynomial {\em reduced modulo} $2$ was obtained by Fuchs and Tabachnikov \cite{FT}. 

\bigskip

\item Inequality c) was proved by Tabachnikov \cite{Ta} using Turaev's state model
for the Kauffman polynomial. 

\bigskip

\item Using another approach, Chmutov, Goryunov and Murakami \cite{CGM} proved inequality b) and Chmutov and Goryunov \cite{CG}
proved inequality c). In \cite{CGM,CG,Ta}, inequalities are stated in the more general context of the contact manifold $ST^* \RM^2$. 
Note also that an analogue of inequality b) for transversal knots
in  $ST^* \RM^2$ is proved independently in  \cite{GH} and in \cite{Ta}.

\bigskip

\item All the results mentioned in this section have a counterpart in Lee Rudolph's theory of
quasi-positive links.
See \cite{Ru1, Ru2, Ru3, Ru4, Ru5, Ru6} for analogues of inequalities a), b) and c).

\bigskip

\item Tanaka \cite{Tan} has shown that inequality $c)$ is a consequence of \cite{Yo},
lemma 1, which itself relies on Turaev's state model for the Kauffman polynomial. 

\end{itemize}

\bigskip

\subsection{Non-Sharpness}

Below in this paper, when it is stated that an inequality of the
form

\bigskip 
 
\centerline{(contact isotopy invariant) $\leq$ (topological invariant)}

\bigskip

is {\em  not sharp}, this means that there exists some topological
knot types $K$ such that the supremum of all the values of
this contact isotopy invariant computed on all Legendrian
representatives of $K$ is less than the value of the
topological invariant computed on $K$.

\bigskip

\begin{itemize}

\item  Using topological methods, Y. Kanda \cite{Ka} has computed
the maximal Bennequin number realizable by a Legendrian representative of some Pretzel knots, showing that the bound 
$tb(l) \leq 2.g_3(l)-1$ is not sharp for these knot types.
The same result follows from Rudolph's \cite{Ru3}, modulo the identification in
\cite{Ru6} of $TB(K)=max\{tb(l) \hskip2mm ; \hskip2mm l 
\textrm{ has topological type} \hskip2mm K\}$ with the invariant $q(K)$ (which was defined, using
another symbol, in \cite{Ru1}). See also \cite{Ru4}.  

\bigskip

\item J. Epstein \cite{Ep} and L. Ng \cite{Ng} have conjectured
  non-sharpness of inequality c) for the knot $8_{19}$. This was
  proved by J. Etnyre and K. Honda \cite{EH}, as a byproduct of 
their classification of Legendrian torus knots.   

\bigskip

\item Sharpness has been established for some specific knot types (see,
for exemple, \cite{Tan,Ep,Ng}).

\bigskip

\item Inequality $a)$ is not optimal already
{\em at the level of concordance classes }  (see \cite{Fe}).

\end{itemize}

\bigskip

\section{Knot polynomials.}
Here, the precise definition of the topological invariants $e_P$ and
$e_Y$ is given. We follow the normalization of \cite{Kau}, pp 215-222.

To any regular oriented knot projection we associate $R$, 
a Laurent polynomial 
in the variables $z$ and $a$ defined by
the following skein relations:
$$R(\rond)=\frac{ a-a^{-1}}{ z}$$
$$R(\fcpo)-R(\fcmo)=z\cdot R(\ouvertVff)$$
$$R(\boucled)=a \cdot R(\droit)$$
$R$ is a regular isotopy invariant and
the HOMFLY polynomial $P(z,a)=a^{-w}R(z,a)$ (where $w=\sharp \fcpo
- \sharp \fcmo$) is a knot invariant.
The least exponent of the variable $a$ in $P$ is
denoted by $e_P$. It is known that $P$ is independent of the orientation, and that $e_P+1$ is an additive
knot invariant with respect to connected sum.

\bigskip

{\bf Example.} $P(\postref)=a^{-3}(\frac{a-a^{-1}}{z})(2a-a^{-1}+az^2)$, hence
$e_P(\postref)=-5$.

\bigskip

 To any regular knot projection we associate $D$, 
a Laurent polynomial of
the variables $z$ and $a$
defined by
the following skein relations:
$$D(\rond)=\frac{ a-a^{-1}}{ z}+1$$
$$D(\fcp)-D(\fcm)=z\cdot(D(\ouvertV) -D(\ouvertH))$$
$$D(\boucled)=a \cdot D(\droit)$$
$D$ is a regular isotopy invariant and
the Kauffman polynomial  $Y(z,a)=a^{-w}D(z,a)$ is a knot invariant.
The least exponent of the variable $a$ in $Y$ is
denoted by $e_Y$. It is known that $e_Y+1$ is an additive knot invariant 
with respect to connected sum.

\bigskip

{\bf Example.} $Y(\postref)=a^{-3}(1+\frac{a-a^{-1}}{z})(2a-a^{-1}+z-a^{-2}z+az^2-a^{-1}z^2)$ 
hence $e_Y(\postref)=-6$.

\bigskip

\section{The Jaeger formula.}
This formula (see \cite{Kau}, pp 219-222) shows that
the Kauffman polynomial of some knot can be computed from the
HOMFLY polynomials of the knots obtained by ``splicing'' a regular
projection of this knot at some crossings.

\bigskip

Consider a link diagram $K$.
A {\em state} $\sigma$ is the following data: A link $K_\sigma$
obtained from $K$
by splicing  some of the crossings (\fcp is modified into \ouvertH or \ouvertV, or is left unchanged), 
and an orientation of $K_\sigma$.
A state $\sigma$ being given, a {\em local weight} is associated the
each crossing $x$ of $K$. If $x$ does not belong to the spliced
crossings, then the
local weight of $x$ is one. 
Consider now an $x$ that belongs to the spliced crossings and suppose that $x=\fcp$ before splicing. 
There are $8$ possible local pictures. If $x$ is spliced to
$\ouvertVn$ then the weight of $x$ is $(t-t^{-1})$. If $x$ is spliced
to $\ouvertHp$, then the weight of $x$ is
$-(t-t^{-1})$. The weight of $x$ vanishes in all remaining possibilities.
The weight of $\sigma$, denoted by $[K,\sigma]$, is the product of all
these local weights.
Denote by $r_\sigma$ the degree of the Gauss map (Whitney index) of
the oriented plane curve underlying the knot diagram $K_\sigma$.

\bigskip

{\bf Theorem.} (Jaeger)
$D(K)(t-t^{-1},a^2t^{-1})= \sum_{\sigma} (ta^{-1})^{r_\sigma}
[K,\sigma]R(K_\sigma)(t-t^{-1},a).$

\bigskip

\section{The Legendrian version of the Jaeger formula.}
It is a reformulation of the formula above in terms of the projection
of Legendrian knots in the plane $(q,u)$, called {\em fronts}.
These are not regular projections.
A generic front has transverse self-intersections \cross and semi-cubic
cusps like  \cuspg or \cuspd. It has no vertical tangent. A typical front is \bifly.
A generic front completely determines the Legendrian knot which lies above, hence {\em
in the sequel a Legendrian knot $l$ is identified with its front when there is no ambiguity.}

\bigskip

To the (generic) front of some Legendrian knot $l$, a generic knot diagram,
called the {\em morsification}
of the front, is associated 
by the following rule:
Each crossing \cross is modified to \fcp.
Each cusp pointing leftward \cuspg is modified to \gauche.
Each cusp pointing rightward \cuspd is modified  to \boucled.
For example, the morsification of \bifly is \trifi.
\bigskip

{\bf Claims.}
The morsification of the front of $l$ is such that:
\begin{itemize}
\item The corresponding knot has the topological type of $l$. 
\item The Whitney index of the morsification is $r=-\mu(l)$.
\item The writhe of the morsification is $w=-tb(l)$.
\item The regular isotopy type of the morsification is invariant under
  Legendrian isotopy. 
\item[$\square$]
\end{itemize}

\bigskip

$R$ and $D$ are 
defined for regular knot diagrams.
Observe that $D$ is defined for unoriented diagrams, and that inverting the orientation leaves $R$,
$tb$ and $w$ invariant (but changes the sign of $\mu$). In the sequel, 
$R(l)$ (resp. $D(l)$) denotes the 
polynomial computed by applying skein relations to the morsification of the front of $l$.
However this is the same as the polynomial computed applying skein relations to the (generically 
regular)  projection of $l$ in the plane $(p,q)$.

\bigskip

{\em  A state} of the front of $l$ consists in the following data: A front $l_\sigma$ 
obtained from the 
one of $l$ by splicing some crossings (\cross can be modified to \ouvertH, or 
 \cuspd  \cuspg, or left unchanged),
and the choice of an orientation of the resulting $l_\sigma$.
A state $\sigma$ of $l$ being given, to each crossing $x=$\cross of $l$, 
a local weight is associated.
If $x$ is left unspliced, then its weight is one.
Suppose now that $x$ belongs to the spliced crossings. There are $8$ possible local
pictures. If $x$ is spliced to \cuspdn  \cuspgp, then its weight is $ta^{-2}(t-t^{-1})$.
If $x$ is spliced to \ouvertHp  then its weight is $t^{-1}-t$.
The weight of $x$ vanishes in all remaining possibilities.
The weight of $\sigma$, denoted by $[l,\sigma]$ is the product of all the local weights.
Denote by $\sharp \cuspgp$ (resp. $ \sharp \cuspdn$) the number of
cusps of $l_\sigma$ which point leftward 
(resp. rightward) and which are oriented upward (resp. downward).
\bigskip

Using this language, the Jaeger formula translates to (see the proof below):

$$ (LJ) \hskip1cm D(l)(t-t^{-1},a^2t^{-1})= \sum_{\sigma} 
(at^{-1})^{\sharp \cuspgp +\sharp \cuspdn}
[l,\sigma]R(l_\sigma)(t-t^{-1},a).$$

\bigskip

{\bf Example 1.} 
$R(\flying)=R(\infi)= \frac{a^2-1} {z}$, and
$D(\flying)=D(\infi)=a+\frac{a^2-1  }{z}$.
There are two states for $\flying$ (the two possible orientations).
Hence the right hand side of $(LJ)$ is:
$$(ta^{-1})^0 \frac{a^2-1}{ (t-t^{-1})} + (at^{-1})^2
\frac{a^2-1} {(t-t^{-1})}.$$ 
This is equal to $D(t-t^{-1}, a^2t^{-1})=
(a^2t^{-1})(1+\frac{a^2t^{-1}-a^{-2}t}{t-t^{-1}})$, as expected.

\bigskip

{\bf Example 2.}
$R(\bifly)=R(\trifi)=\frac{a^3-a} {z}$.
$D(\bifly)=D(\trifi)=a^2+\frac{a^3-a  }{z}$.
There are 4 states whose weights do not vanish: \statun, \statdeux, \stattrois \stattrois,
and \stattrois.
The right hand side of $(LJ)$ is:
$$(at^{-1}) \frac{a^3-a} {t-t^{-1}}
+(at^{-1})  \frac{a^3-a} {t-t^{-1}}
+(at^{-1})^4  (ta^{-2})(t-t^{-1}) (\frac{a^2-1} {t-t^{-1}})^{2}
+(at^{-1})^2 (t^{-1}-t) \frac{a^2-1} {t-t^{-1}}.$$
This is equal to $D(t-t^{-1},a^2t^{-1})=(a^2t^{-1})^2(1+\frac{a^2t^{-1}-a^{-2}t}{t-t^{-1}})$, 
as expected.

\bigskip

{\it Proof of (LJ).} Consider some Legendrian knot $l$, and the knot diagram $K$ obtained by 
rounding all the cusps of the front of $l$ (\cuspg becomes \gauche and \cuspd becomes \droit).
Denote by $\nu$ half the number of cusps of $l$. By the axioms for $D$, 
$$D(l)(t-t^{-1},a^2t^{-1})=(a^2t^{-1})^{\nu}D(K)(t-t^{-1},a^2t^{-1}).$$
There is a one-to-one correspondence between the states of $l$ and the states of $K$. 
Writing the Jaeger formula for $K$ in terms of $l$ will give $(LJ)$:
Consider some state $\sigma$ of $l$. Denote by $\nu_\sigma$ half the
number of cusps of $l_\sigma$ and by $V$ (resp. by $H$) the number of crossings
of $l$ (or of $K$) that are spliced vertically (resp. horizontally) in $\sigma$.
The following relations hold: $R(K_\sigma)=a^{-\nu_\sigma}R(l_\sigma)$, 
$[K,\sigma]=(-1)^H(t-t^{-1})^{V+H}$, $\nu=\nu_\sigma-V$, and 
$\nu_\sigma-r(K_\sigma)={\sharp \cuspgp +\sharp \cuspdn}$.
Plug this into the expression of $D(l)$ above:
$$ D(l)(t-t^{-1},a^2t^{-1})= \sum_{\sigma} 
(at^{-1})^{\sharp \cuspgp +\sharp \cuspdn}
(-1)^H(a^2t^{-1})^{-V}(t-t^{-1})^{V+H}R(l_\sigma)(t-t^{-1},a). $$
This is $(LJ)$.
$\square$
\bigskip

\section{Inequality $c)$ follows from inequality $b)$}
\label{b_implique_c}

Since $tb=-w$, inequality $b)$
is equivalent to the fact that there is no negative 
power of $a$ occurring in $a^{-|\mu|}R(l)$, i.e., {\em it is a genuine polynomial in $a$}.
Similarly, we want to prove that $D(l)$ is a genuine polynomial in $a$. 
This is a consequence of the following lemma.

\bigskip

{\bf  Lemma.} The contribution of each state in the right hand side of $(LJ)$ is a genuine 
polynomial  in $a$.

\bigskip

{\it Proof.} 
Consider a state $\sigma$ of $l$.
Denote by $V$ the number of crossings that are spliced to \cuspdn \cuspgp, and by
$H$ the number of crossings that are spliced to \ouvertHp.
The contribution of $\sigma$ is 
$$(at^{-1})^{\sharp \cuspgp + \sharp \cuspdn}
(ta^{-2})^{V}(-1)^{H}(t-t^{-1})^{V+H}R(l_\sigma)(t-t^{-1},a).$$
The least exponent of $a$ in $R(l_\sigma)(t-t^{-1},a)$ is not less than 
$|\mu(l_\sigma)|$, by inequality $b)$.
Denote by $E$ the least exponent of $a$ in the contribution of $\sigma$. 
$E \geq \sharp \cuspgp  + \sharp \cuspdn -2\cdot V +|\mu|$.
On the other hand $\mu=\sharp \cuspgp  - \sharp \cuspdn$, hence $E \geq 2 
(\sharp \cuspgp - V) +|\mu|-\mu$.
Since splicing \cross to \cuspdn \cuspgp creates one \cuspgp, $V$ is not bigger than 
$\sharp \cuspgp$, and hence $E\geq0$. $\square$

\bigskip

{\bf Remark about this proof.} As explained in \cite{FT}, inequality b) has a ``simple'' and natural 
proof by \cite{Be} and \cite{Mo} or \cite{FW}, much simpler than the known proofs of a) for instance. 
Since the Jaeger formula is also proved (\cite{Kau}) by "elementary" means (like checking its invariance 
under the Reidemeister moves), this proof of $c)$ is, in my opinion, simple and natural.
I find it remarkable that the Jaeger formula fits so well between b) and c).
\bigskip
 
{\bf Remark.} Like $a)$, $b)$ follows from a more general inequality about transverse knots
(see \cite{Be,Ta,GH}).
The proof above, which lacks of a natural transverse counterpart, 
seems to indicate that $c)$ is an inequality about Legendrian knots only. 

\bigskip

\section{Relationship between $g_4$ and $e_P$.}
It is proved here that inequality $a)$ can be stronger than
inequality $b)$.

\bigskip

{\bf Proposition.} 
The difference between $e_P$ and $2 \cdot g_4 -1$ can be arbitrarily negative or positive.

\bigskip

{\bf Corollary.} Inequality b) is not
sharp, i.e. $$max\{tb(l)+ |\mu(l)| \hskip2mm ; \hskip2mm l 
\textrm{ has topological type} \hskip2mm K\}<e_P(K)$$ for some knot types $K$.

\bigskip

{\bf Remark.} $2 \cdot g_4 -1 < e_P$ seems much more difficult to realize than the converse:
The tables indicate no
contradiction to 
$e_P(K) \leq 2 \cdot g_4(K)-1$ for the $84$ first knots (which arise from diagrams with
less than $10$ crossings). 

\bigskip

{\bf Question.} This leaves the question of Morton \cite{Mo} open: Is it
true that $e_P(K) \leq 2 \cdot g_3(K)-1$? (Recall that $g_3$ denotes
the genus). 
This inequality is true for alternating
knots and for positive knots, as proved in \cite{Cr}, and for knots
with braid index less than 4 \cite{DM}. It was checked by Alexander
Stoimenow for all knots which admit a diagram with less than 17 crossings.

\bigskip

{\it Proof of the proposition.}
Consider some knot $K$ such that $e_P+1$ is negative.
For instance $K=\postref$.
Since $e_P+1$ is additive under connected sum, there exist knots with arbitrarily negative $e_P$.
On the other hand, $g_4$ is never negative. Hence $e_p-(2 \cdot g_4-1)$ can be 
made arbitrarily negative.

\bigskip

Let K be the closure of the braid 
$\sigma$= \blanc \sbb \sbb \scc \scc \saa \saa \sbb \sccc \sa \sa \sa \blanc.
This knot admits a projection with ten crossings.
Changing the first $\sbb$ \blanc of $\sigma$ to $\sbbb$ \blanc, one gets 
a braid whose closure is the trivial knot. Hence $g_4(K) \leq 1$ (it is in fact $1$).
Computation shows that $e_P(K)=3$ (hence $K$ is an example for which $b)$ is not sharp).
Denote by $K^{\sharp d}$ the connected sum of $k$ copies of $K$. $g_4(K^{\sharp d}) \leq d$, 
and $e_P(K^{\sharp d})=d(3+1)-1$. Hence $e_P-(2 \cdot g_4 - 1)$ can be made 
arbitrarily positive.  $\square$

\section{Relationships between inequalities $b)$ and $c)$.} \label{ineq}

By section \ref{b_implique_c}, inequality $b)$ implies inequality $c)$. 
However, inequalities $b)$ and $c)$ are independent in the following sense:
$b)$ implies that $tb(l)\leq e_P(l)$. 
Looking in the tables seems to indicate that $e_Y \leq e_P$, and hence that  $tb(l)\leq e_P(l)$
is weaker than $tb(l) \leq e_Y(l)$ (inequality $c)$). This is however not true. 
Among all the
prime knots which admit a diagram with less than 15 crossings (there
are grosso-modo 60.000 of them), there are 22 knots verifying
$e_P < e_Y$. One of the two examples with 12 crossings is the closure 
of the following braid:

\centerline{\Scc \Sbb \Saa \Sbbb \Scc \Sddd \Saaa \Sbb \Scc \Sddd \Sccc \Sbbb \Scc \Saaa 
\Sbb \Saaa \Sdd \Scc 
\Sbb \Saa }

\bigskip

{\bf Corollary.} Inequality $c)$ is not sharp.

\bigskip

{\bf Question.} Is it true that $e_Y \leq e_P$ for alternating knots 
(none of the 13 examples cited above is alternating)?

\bigskip

\section{Acknowledgements}
I learnt the Jaeger formula from Christan Blanchet's lectures at the 
summer school on finite type invariants of knots and
three manifolds, Grenoble, 1999.
All computations of knot polynomials used throughout this paper are due to 
Jim Hoste and Morwen Thistlethwaite, via their  {\em Knotscape} program
\cite{HT}.
Further computer aided example search was programmed by Xavier
Dousson and Alexander Stoimenow.
I thank them all very much.

\end{document}